\documentclass{amsart}
\usepackage{amssymb,latexsym,epic} 
\newlength{\hs} 
\newlength{\agk}
\newlength{\agku}
\newlength{\jmr}
\newlength{\bernd}
\newcommand{\tor}{\mathrm{Tor}}

\newcommand{\cD}{\mathcal{D}}

\newcommand{\cO}{\mathcal{O}}

\newcommand{\cP}{\mathcal{P}}

\newcommand{\supp}{\mathrm{Supp}}

\newcommand{\R}{\mathbb{R}}
\newcommand{\C}{\mathbb{C}}

\newcommand{\Cn}{\C^n}

 \newcommand{\Csn}{{(\C^*)}^n}
\newcommand{\Ksn}{{(K^*)}^n}
 
\newcommand{\cM}{\mathcal{M}}

\begin{document}

\title[Extensions for ``A Convex Geometric Approach...'']{Extensions and  
Corrections for: ``A Convex Geometric Approach to Counting the Roots of a 
Polynomial System$^1$''}
\footnotetext[1]{The paper we are extending and correcting originally appeared in Theoretical Computer Science 133 (1994), pp.\ 105-140,
Elsevier.}

\author{J. Maurice Rojas}
\thanks{This research was completed at MSRI and partially funded by a
an NSF Mathematical Sciences Postdoctoral Research Fellowship and NSF grant 
DMS-9022140.}

\address{Massachusetts Institute of Technology\\
Mathematics Department\\
77 Mass.\ Ave.\\
Cambridge, MA \ 02139, U.S.A. }

\email{rojas@math.mit.edu} 

\date{October 30, 1996} 

\dedicatory{This paper is dedicated to the memory of Laxmi Patel.}

\begin{abstract} 
This brief note corrects some errors in the paper quoted in the title, 
highlights a combinatorial result which may have been overlooked, and points 
to further improvements in recent literature. 
\end{abstract}

\maketitle

\section{Introduction}
\label{sec:intro}
This brief note contains important additional information relating to 
an earlier paper of the present author: ``A Convex Geometric Approach to 
Counting the Roots of a Polynomial System,'' Theoretical Computer Science 133 
(1994), pp.\ 105-140 (henceforth referred to as \cite{convexapp}). 
This paper gave various extensions of the seminal works 
\cite{kus75,bernie,kus76,khocompat} relating root counting for 
systems of polynomial equations to volumes of polyhedra. We will 
report (briefly) on the status of some of these extensions, 
and correct some errors appearing in \cite{convexapp}. 

For the sake of brevity, we will not review any notation or definitions, 
since they are already amply covered in \cite{convexapp} and the 
more recent \cite{toricint}. However, we point out that the latter 
paper is readily available on-line at {\tt http://www-math.mit.edu/\~{}rojas}.

We begin, in the following section, by pointing out a combinatorial result 
from \cite{convexapp} which seems to have gone ignored. (In particular, 
a special case of \cite[Corollary 3]{convexapp} was the main 
result of a paper completed in 1996 by another author!) In the 
next section we then discuss certain problems within root counting, for 
$n\times n$ polynomial systems, which are close (or not so close) to a 
satisfactory solution. We then provide a list of 
corrigenda for \cite{convexapp} in the final section. 

\section{Filling and Counting} 
A combinatorial corollary of the results of \cite{convexapp} gave the 
first known constructive solution of the {\em filling}\/ problem for rational 
polytopes \cite{convexapp,drs}. (This also resolves (for rational 
polytopes) a conjecture of Rolf Schneider on the mixed area measure
\cite{schneider}.) More explicitly, a combinatorial answer is given to the 
following question: 
Given an $n$-tuple of rational polytopes with positive mixed volume, 
which {\em sub}-$n$-tuples of rational polytopes have the same mixed 
volume? A partial answer is contained in 
Corollary 4 of \cite[Page 115]{convexapp} and a full answer 
appears in Corollary 9 on Page 136. The proof comes down to a 
carefully tailored application of Bernshtein's Theorem \cite{bernie}, 
contained in Lemmata 2--4 and Corollary 5 of \cite[Pages 121--123]{convexapp}. 
In particular, this convex geometric problem was first solved by an algebraic 
geometric result. 

More importantly, the filling problem was shown to be equivalent to 
the {\em $\Ksn$-counting}\/ problem. The latter problem was defined in 
\cite{convexapp} as the classification of all subsets 
of coefficients of a given polynomial system (with {\em fixed}\/ monomial term 
structure) whose genericity guarantees that the mixed volume bound is an exact 
root count. (The notion of $\Ksn$-counting was referred to as 
{\em counting}\/ in \cite{convexapp}.) 
The equivalence of filling and $\Ksn$-counting is the content of Lemmata 2 and 
3 of \cite{convexapp}. 

Unfortunately, filling and counting were extended in an inelegant 
way in \cite{convexapp}. This may have been the cause for the 
obscurement of filling and counting. For instance, a friend of the 
author's (who will remain unnamed) wrote an entire paper based 
on a special case of Corollary 3 of \cite{convexapp}. It is 
the author's firm advice to ignore $r$-counting and $(r,s,n)$-filling 
(which are admittedly quite abstruse), and instead follow the 
improved construction of $\mathbf{W}${\bf -counting} of 
\cite{rojaswang,toricint}. 

It should also be remarked that an earlier solution to the $\Csn$-counting 
problem was incomplete: In \cite{johnandme}, it was falsely asserted 
that the {\em Vertex Coefficient Theorem}\/ gave a complete solution. This was 
later corrected in the latter author's M.S. thesis \cite{myms} and the complete 
solution seems to have first appeared in \cite[Lemmata 2 and 3]{convexapp}. 
Finally, we remark that $\Ksn$-counting is also a much simpler criterion than 
the {\em ID cover}\/ of \cite{johnandme}. 

\section{Extensions --- Complete and Incomplete} 
Recall that $\Csn\!:=\!(\C\!\setminus\!\{0\})^n$ is sometimes 
referred to as the (complex) {\em algebraic torus}. 
The BKK bound \cite{kus75,bernie,kus76,khocompat} was a 
beautiful result discovered almost two decades ago in a seminar 
of V. I. Arnold. This result gave an upper bound on the number 
of isolated roots in the algebraic torus (of $n$ polynomial equations in $n$ 
unknowns) in terms of volumes of $n$-dimensional polyhedra. 

These upper bounds also possesed an extremely 
important property: they were the best possible, given only the monomial 
term structure.  In other words, if one fixed which monomial terms 
appeared in a polynomial system, the resulting convex geometric formula would 
fail to be an {\em exact}\/ root count only on a positive codimension algebraic 
subset of the coefficient space. (The terminology ``generically exact'' 
is also sometimes used in this respect.) Aside from a result of 
F.\ Minding for the case of two equations in two unknowns \cite{min}, 
this sort of optimality had {\em never}\/ been attained by any 
previous upper bound.

Four natural extensions of the preceding result immediately come to mind:
\begin{enumerate}
\item{Counting the {\em exact} number of roots (and not just a tight upper 
bound).}
\item{Extending to algebraically closed fields other than $\C$ (and 
in particular, to positive characteristic).}
\item{Counting the number of roots in $\Cn$ and in subregions other than 
$\Csn$.}
\item{Getting information about the higher (co)homology structure 
of other locally complete intersections.}
\end{enumerate}

We now briefly point out what has been done from 1994 to 1996 for these 
problems:
\begin{enumerate}

\item{We first remark that there are now precise combinatorial and 
algebraic conditions for when the BKK bound fails to be an exact root 
count. Combinatorial conditions first appeared in \cite{johnandme} 
and were then further refined in \cite{convexapp,rojaswang,toricint}; 
algebraic conditions first appeared in \cite{bernie} and were then refined in 
terms of the sparse resultant in \cite{polyhomo} (as well as in independent  
work of the present author). The latter refinement was then extended 
further (and corrected) in \cite{toricint}. 

There are, of course, many ways to count the exact number of roots 
directly with commutative algebra and Gr\"obner bases. For instance, 
a particularly nice approach (which works over $\R$ as well) 
is a recent extension \cite{prs93} of {\em Hermite's}\/ method 
\cite{meat}. 

However, there are more recent methods based on  
toric geometry which make more refined use of the monomial 
term structure of a given problem. For example, \cite{lamina} gives a new  
method, based on the sparse resultant, to count the exact number of roots in 
the algebraic torus. It is interesting to note that the original 
BKK bound (at worst) fails to be an exact root count on a codimension 1 
subset of the space of coefficients. The method in \cite{lamina} fails only 
on a codimension $\geq\!2$ subset, and an extension which {\em always}\/ 
works has just been completed \cite{twist}. Thus, there are now convex 
geometric methods (augmented by the sparse resultant) to count the exact 
number of roots in $\Csn$. It is hoped that these methods will prove 
significantly faster for exact root counting than current 
Gr\"obner basis methods. } 

\item{Convex geometric root counting can be done over {\em any}\/ algebraically
closed field --- not just $\C$. This began with Danilov's more abstract 
framework \cite{dannie} for the BKK bound, and was further refined in 
\cite{convexapp} and \cite{toricint}. The combinatorial and 
algebraic conditions for exactness of the BKK bound also hold over any 
algebraically closed field 
\cite{toricint}. Finally, the aforementioned extensions of optimal upper 
bounds to exact root counts (via sparse resultants) work over any 
algebraically closed field as well.} 

\item{It is indeed possible to get optimal convex geometric upper bounds on
the number of roots in all of $\Cn$. This was first considered in 
\cite{kho78}, for certain polynomial systems. Suboptimal upper bounds, valid
for all polynomial systems, were then 
derived in \cite{convexapp} and \cite{rojaswang,liwang}. The last 
two papers gave, respectively, combinatorial and complex geometric 
conditions for when their convex geometric bounds were optimal. (It 
should also be noted that \cite{rojaswang} gave tight upper bounds on the 
number of roots in affine space minus an arbitrary union of coordinate 
hyperplanes, over any algebraically closed field.) The first optimal 
bounds for $\Cn$ minus an arbitrary union of coordinate hyperplanes, 
holding for {\em all}\/ polynomial systems, appeared in \cite{hsaff}. In fact, 
their results held in greater generality: Any Boolean vanishing condition on 
the coordinates $x\!=\!(x_1,\ldots,x_n)$, e.g., 
$(x_1\!=\!0)\wedge(x_3\!\neq\!0)\vee\cdots$, was allowed and 
such roots could also be (generically) counted convex geometrically. These  
results 
were then extended to arbitary algebraically closed fields, and an alternative 
formula derived, in \cite{toricint}. The preceding algebraic and 
combinatorial conditions for exactness were also extended to the 
case of affine space minus an arbitrary union of coordinate hyperplanes 
in the same paper. }
\item{Convex geometric upper bounds on the degree of certain 
positive-dimensional varieties were derived in \cite{convexapp}. These 
results overlapped slightly with the deeper results of \cite{dankho} 
on finding the {\em mixed Hodge structure}\/ of a variety via convex 
geometry. For example, a special case of the latter work gave a 
convex geometric formula for the Euler characteristic of certain 
(generic) subvarieties of $\Csn$. Another often overlooked example is their  
(generically valid) computation of arithmetic genus via the number of lattice 
points in a polyhedron. Combined with \cite{toricint}, 
it now appears that the results of \cite{dankho} can be {\em completely}\/ 
extended to affine space and arbitrary algebraically closed fields. However, 
the question of finding {\em precise}\/ algebraic or combinatorial conditions 
for when their more general formulae hold is still open. } 
\end{enumerate}

In closing, we add that there are still (as of 1996) no {\em proven}\/ convex 
geometric formulae for the maximal number of {\bf real} roots. An important 
step toward this goal is the conjectural formula of Sturmfels 
\cite{realstu91}, later simplified by Itenberg and Roy \cite{ilyamf}, 
which attempts to generalize Descartes' rule to higher dimensions. 
An interesting explicit formula for the expected number of real 
roots of certain {\em random}\/ sparse polynomial systems appears in 
\cite{myaverage}. 

\section{Corrections to \cite{convexapp}}
\begin{itemize}
\item[A Frequent Typo:]{In any places, nearby semicolons and commas should be 
reversed to correct the appearance of $\cM(P,n)$, $\cM(\cP;\Delta_s,n-k)$, and 
$\cM_r(\cP;\Delta_s,n-k)$. }
\item[Page 117, Line -16:]{In general, the definition given for 
intersection multiplicity is only an {\em upper bound}. However, 
when $W$ is a proper $0$-dimensional component, the formula is exact 
and can be further simplified to the {\em dimension}\/ 
(as a $K$-vector space) of the same $R$-module. To correctly define 
intersection multiplicity for a proper {\em positive}-dimensional 
component, it is necessary to use $\tor(\cdot)$ as in Serre's 
construction \cite{ifulton2}. }
\item[Page 118, Line 2:]{It should have been mentioned 
that throughout the paper, $\deg W$ actually means the 
degree of the {\em reduced}\/ variety defined by $W$.}
\item[Page 115, Line -1:]{The second to last sentence should 
end with ``...in $\Csn$.''.}
\item[Page 119, Line 15:]{Fact (2) is incorrect. The proper 
statement involves a related (canonically defined) intersection 
of toric divisors and appears in \cite[Corollary 2]{toricint}. }
\item[Page 119, Line 18:]{``...$\supp(\cD_i)$ contains...'', not 
``...supported precisely on...''}
\item[Page 119, Line 24:]{The line bundles $\cO(\cD_i)$ are 
{\em not}\/ ample in general. (I thank Professor William Fulton for 
pointing this out to me.) However, through an algebraic homotopy argument, 
one can still obtain the ``numerical ampleness'' assertion of the final 
sentence of the proof. This is done in detail in the proof of Theorem 3 of 
\cite{toricint}. }
\item[Section 2.6]{The results of this entire section are considerably 
simplified and improved in \cite{rojaswang,hsaff,toricint}.}
\item[Page 138, Line -15:]{``...when one goes...'', not ``when goes''. }
\item[Reference \text{[13]}:]{...is now entitled ``How to Fill a Mixed 
Volume.''}
\end{itemize}

\bibliographystyle{amsalpha}

\begin{thebibliography}{A}

\bibitem[Ber75]{bernie} Bernshtein, D. N., {\it ``The Number of Roots of a
System of Equations,"} Functional Analysis and its Applications (translated
from Russian), Vol. 9, No. 2, (1975), pp.\ 183--185.

\bibitem[CR91]{johnandme} Canny, John F. and Rojas, J. Maurice,
{\it ``An Optimal Condition for Determining the Exact Number of Roots
of a Polynomial System,"} Proceedings of ISSAC '91 (Bonn, Germany), ACM
Press (1991), pp.\ 96--102.

\bibitem[Dan78]{dannie} Danilov, V. I., {\it ``The Geometry of Toric
Varieties,"} Russian Mathematical Surveys, 33 (2), pp.\ 97--154, 1978.

\bibitem[DK87]{dankho} Danilov, V. I., and Khovanskii, A., {\it ``Newton
Polyhedra and an Algorithm for Computing Hodge-Deligne Numbers,"} Math.
USSR Ivzestiya, Vol.\ 29 (1987), No.\ 2.

\bibitem[DRS96]{drs} Dalbec, J., Rojas, J. M., and Sturmfels, B., {\it
``How to Fill a Mixed Volume,"} manuscript, Massachusetts Institute of
Technology, 1996.

\bibitem[Ful84]{ifulton2} Fulton, William, {\it Intersection Theory,} 
Springer-Verlag, 1984.

\bibitem[Her56]{meat} Hermite, C., {\it ``Sur le Nombre des Racines d'une 
\'Equation Alg\'ebrique Comprise Entre des Limites Donn\'ees,''} 
J. Reine Angew. Math., vol.\ 52, pp.\ 39--51 (1856). 

\bibitem[HS95]{polyhomo} Huber, Birkett and Sturmfels, Bernd, {\it ``A
Polyhedral Method for Solving Sparse Polynomial Systems,"} Mathematics of
Computation, 64, pp.\ 1541--1555, 1995.

\bibitem[HS96]{hsaff} \underline{\hspace{\hs}}, {\it ``Bernshtein's Theorem in 
Affine Space,"} Discrete and Computational Geometry, to appear, 1996.

\bibitem[IR95]{ilyamf} Itenberg, Ilya and Roy, Marie-Francoise,
{\it ``Multivariate Descartes' Rule,"} Preprint (1995), IRMAR,
Universit\'e de Rennes I, France.

\bibitem[Kho77]{khocompat} Khovanskii, Askold G., {\it ``Newton Polyhedra and
Toroidal Varieties,"} Functional Anal. Appl., 11 (1977), pp.\ 289--296.

\bibitem[Kho78]{kho78} \underline{\hspace{\agk}}, {\it ``Newton Polyhedra and 
the Genus of Complete Intersections,"} Functional Analysis  
(translated from Russian), Vol. 12, No. 1, January--March (1978), pp.\ 51--61.

\bibitem[Kus75]{kus75} Kushnirenko, A. G., {\it ``A Newton Polytope and the
Number of Solutions of a System of k Equations in k Unknowns,"} Usp. Matem.
Nauk., 30, No.\ 2, pp.\ 266--267 (1975).

\bibitem[Kus76]{kus76} \underline{\hspace{\agku}}, {\it ``Newton Polytopes and 
the B\'ezout Theorem,"} Functional Analysis and its Applications (translated 
from Russian), vol.\ 10, no.\ 3, July--September (1976), pp.\ 82--83.

\bibitem[LW96]{liwang} Li, T. Y. and Wang, Xiaoshen, {\it ``The BKK Root
Count in $\Cn$,"} Mathematics of Computation, October, 1996. 

\bibitem[Min41]{min} Minding, F., {\it ``\"Uber die Bestimmung des
Grades der durch Elimination hervorgehenden Gleichung,"} J. Reine und Angew.
Math. (1841), 22, pp.\ 178--183.

\bibitem[PRS93]{prs93} Pedersen, P., Roy, M.-F., and Szpirglas, A.,
{\it ``Counting Real Roots in the Multivariate Case,"}
Computational Algebraic Geometry, edited by Eyssette and Galligo, Progress
in Mathematics 109, pp.\ 203--224, Birkhauser, 1993.

\bibitem[Roj91]{myms} Rojas, J. Maurice, {\it ``An Optimal Condition for
Determining the Exact Number of Roots of a Polynomial System,"} M.S.
thesis, C.S. Division, U. C. Berkeley, 1991.

\bibitem[Roj94]{convexapp} \underline{\hspace{\jmr}}, {\it ``A Convex Geometric
Approach to Counting the Roots of a Polynomial 
System,"} Theoretical Computer Science (1994), vol.\ 133 (1), pp.\ 105--140.

\bibitem[Roj96a]{myaverage} \underline{\hspace{\jmr}}, {\it ``On the Average 
Number of Real Roots of Certain Random Sparse,''} pp.\ 689--699, {\it The
Mathematics of Numerical
Analysis}, Lectures in Applied Mathematics, vol.\ 32 (1996), edited by Jim
Renegar, Mike Shub, and Steve Smale, American Mathematical Society. (Also
available on-line at {\tt http://www-math.mit.edu/\~{}rojas}.)

\bibitem[Roj96b]{toricint} \underline{\hspace{\jmr}}, {\it ``Toric  
Intersection Theory for Affine Root Counting,''} submitted to the  
Journal of Pure and Applied Algebra. (Also available on-line at 
{\tt http://www-math.mit.edu/\~{}rojas}.)

\bibitem[Roj96c]{lamina} \underline{\hspace{\jmr}}, {\it ``Toric Laminations,   
Sparse Generalized Characteristic Polynomials,  
and a Refinement of Hilbert's Tenth Problem,''} Proceedings of the Rio de 
Janeiro Foundations of Computational Mathematics Conference (January 1997), 
Springer-Verlag, to appear. 

\bibitem[Roj96d]{twist} \underline{\hspace{\jmr}}, {\it ``Polynomial 
Equation Solving and Twisted Chow Forms,''} manuscript, MIT. 

\bibitem[RW96]{rojaswang} Rojas, J. M., and Wang, Xiaoshen, {\it ``Counting
Affine Roots of Polynomial Systems Via Pointed Newton Polytopes,"} Journal
of Complexity, vol.\ 12, June (1996), pp.\ 116--133. 

\bibitem[Sch94]{schneider} Schneider, Rolf, {\it Convex Bodies: The
Brunn-Minkowski Theory,} Encyclopedia of Mathematics and its Applications,
v.\ 44, Cambridge University Press, 1994.

\bibitem[Stu91]{realstu91} Sturmfels, Bernd,
{\it ``On the Number of Real Roots of a Sparse Polynomial System,"}
Hamiltonian and Gradient Flows: Algorithms and Control, (ed.\ A. Bloch),
Fields Institute Communications, Vol. 3, American Math. Soc., Providence,
RI,
1991, pp.\ 137--143.


\end{thebibliography}

\end{document}